\documentclass{amsart}
\usepackage{amsmath, amssymb, latexsym, setspace}
\title{Flag f-vectors of colored complexes}
\author{Andrew Frohmader}
\address{Department of Mathematics, 581 Malott Hall, Cornell University, Ithaca, NY 14853-4201}
\email{froh@math.cornell.edu}

\newtheorem{theorem}{Theorem}

\newtheorem{definition}[theorem]{Definition}

\def\proof{\smallskip\noindent {\it Proof: \ }}
\def\endproof{\hfill\ensuremath{\square}\medskip}

\begin{document}

\maketitle

\begin{abstract}
It is shown that conditions stronger in a certain sense than color-shifting cannot be placed on the class of colored complexes without changing the characterization of the flag f-vectors.
\end{abstract}

In the late 1970s, Stanley \cite{stanley} showed that two particular classes of simplicial complexes have equivalent characterizations of their flag f-vectors.  Several years later, Bj\"{o}rner, Frankl, and Stanley \cite{bfs} showed that two additional classes of simplicial complexes shared this same characterization.  Unfortunately, no one has a characterization for any of these classes of simplicial complexes, but we only know that characterizing one would suffice for all four.

The two additional classes of simplicial complexes included in the equivalence of Bj\"{o}rner, Frankl, and Stanley are each proper subsets of one of the classes of complexes in Stanley's original paper.  Thus, the paper of Bj\"{o}rner, Frankl, and Stanley could be thought of as progress toward a solution by narrowing the class of complexes to consider.  In this paper, we show in Theorem~\ref{nolocal} that extending this approach to a solution of the problem by further narrowing one of the classes of complexes in a certain sense is impossible.

Recall that a \textit{simplicial complex} $\Delta$ on a vertex set $W$ is a collection of subsets of $W$ such that (i) for every $v \in W$, $\{v\} \in \Delta$ and (ii) for every $B \in \Delta$, if $A \subset B$, then $A \in \Delta$.  The elements of $\Delta$ are called \textit{faces}.  A face on $i$ vertices is said to have \textit{dimension} $i-1$, while the dimension of a complex is the maximum dimension of a face of the complex.

The \textit{$i$-th f-number} of a simplicial complex $\Delta$, $f_{i-1}(\Delta)$ is the number of faces of $\Delta$ on $i$ vertices.  The \textit{f-vector} of $\Delta$ lists the f-numbers of $\Delta$.  One interesting question to ask is which integer vectors can arise as f-vectors of simplicial complexes.  Much work has been done toward answering this for various classes of simplicial complexes.  For example, the Kruskal-Katona theorem \cite{kruskal, katona} characterizes the f-vectors of all simplicial complexes.

In this paper, we wish to deal with colored complexes, where the coloring provides additional data.  A \textit{coloring} of a simplicial complex is a labeling of the vertices of the complex with colors such that no two vertices in the same face are the same color.  Because any two vertices in a face are connected by an edge, this is equivalent to merely requiring that any two adjacent vertices be assigned different colors.  If the set of colors has $n$ colors, we refer to the colors as $1, 2, \dots , n$.  The set of colors is denoted by $[n] = \{1, 2, \dots , n \}$.  The \textit{color set} of a face is the subset of $[n]$ consisting of the colors of the vertices of the face.  The Frankl-F\"{u}redi-Kalai \cite{ffk} theorem characterizes the f-vectors of all simplicial complexes that can be colored with $n$ colors.

We wish to use a refinement of the usual notion of f-vectors. The \textit{flag f-numbers} of a colored simplicial complex $\Delta$ on a color set $[n]$ are defined by, for any subset $S \subset [n]$, $f_S(\Delta)$ is the number of faces of $\Delta$ whose color set is $S$.  The \textit{flag f-vector} of $\Delta$, $f(\Delta)$, is the collection of the flag f-numbers of $\Delta$ for all subsets $S \subset [n]$.

This is a refinement of the usual notion of f-numbers and the f-vector of a complex.  The relation between the f-numbers and the flag f-numbers is that the former ignores the colors, and can be computed from the latter as
$$f_{i-1}(\Delta) = \sum_{|S| = i} f_S(\Delta).$$

One can ask which nonnegative integer vectors can arise as the flag f-vectors of colored simplicial complexes.  It can help to define the flag h-vector of a complex by
$$h_S(\Delta) = \sum_{T \subset S} f_T(\Delta) (-1)^{|S| - |T|}.$$
The flag h-vector of a complex contains the same information as the flag f-vector, and is easier to work with in some cases.  If given the flag h-vector, we can recover the flag f-vector by
$$f_S(\Delta) = \sum_{T \subset S} h_T(\Delta).$$

One approach is to try to find some bounds.  Walker \cite{walker} showed that the only linear inequalities on the flag f-numbers of simplicial complexes are the trivial ones, namely, that all flag f-numbers are non-negative.  He also computed all linear inequalities on the logarithms of the flag f-numbers of a simplicial complex.  These give inequalities on the products of flag f-numbers.  Walker's bounds are not sharp, but they are enough to settle the case of two colors.  A proposed nonnegative integer flag f-vector corresponds to a non-empty two-colored simplical complex if and only if $f_{\emptyset}(\Delta) = 1$ and $f_1(\Delta)f_2(\Delta) \geq f_{12}(\Delta)$.  This does not settle the problem for more colors, however.

Another paper by the author \cite{threecolor} characterizes the flag f-vectors of three-colored complexes.  The characterization is somewhat complicated, and essentially consists of trying to construct a complex with the given flag f-vector in several ways, with typically around six or ten things to try.  If any of them work, then it gives a colored complex with the desired flag f-vector.  If none of them work, then the paper shows that there is no colored complex with the desired flag f-vector.

Stanley \cite{stanley} showed that the flag h-vector of a balanced Cohen-Macaulay complex is the flag f-vector of a simplicial complex and vice versa.  That is, if a Cohen-Macaulay complex has dimension $n-1$ and can be colored with $n$ colors, then there is a simplicial complex that can be colored with $n$ colors whose flag f-vector is the flag h-vector of the Cohen-Macaulay complex.

Bj\"{o}rner, Frankl, and Stanley \cite{bfs} were able to further restrict both of these classes of complexes.  We need a bit of notation to state the relevant portion of their results.  We can place an arbitrary order on the vertices of each color.  We label the $j$-th vertex of color $i$ as $v_j^i$, so that the vertices of color $i$ are $v_1^i, v_2^i, \dots, v_{f_i(\Delta)}^i$.

\begin{definition}
\textup{Let $\Delta$ be an $n$-colored simplicial complex.  We say that $\Delta$ is \textit{color-shifted} if, for all $b_1 \leq a_1, b_2 \leq a_2, \dots , b_j \leq a_j$, $\{v_{a_1}^{i_1}, v_{a_2}^{i_2}, \dots v_{a_j}^{i_j}\} \in \Delta$ implies $\{v_{b_1}^{i_1}, v_{b_2}^{i_2}, \dots v_{b_j}^{i_j}\} \in \Delta$.}
\end{definition}

\begin{theorem} \cite[Theorem 1]{bfs} \label{colorshift}
Let $\Delta$ be an $n$-colored simplicial complex.  Then there is an $n$-colored, color-shifted simplicial complex $\Gamma$ such that $f_S(\Delta) = f_S(\Gamma)$ for all $S \subset [n]$.
\end{theorem}

Bj\"{o}rner, Frankl, and Stanley called this concept ``compressed" rather than color-shifted.  Furthermore, their proof allowed for a more general notion of coloring where, for example, one could have three colors, but allow a face to have up to 3 vertices of color 1, up to 5 vertices of color 2, and up to 2 vertices of color 3.  In this paper, we focus on the case where only one vertex of each color is allowed in a face.

The converse of Theorem~\ref{colorshift} is trivial, as every color-shifted colored complex is, in particular, a colored complex.  What this does is to say that rather than considering all colored complexes, it suffices to consider only those that are color-shifted.  Likewise, the same theorem of theirs said that rather than considering all Cohen-Macaulay balanced complexes, it suffices to consider shellable balanced complexes.  The problem here is that while four different classes of complexes have equivalent characterizations, none of them have a known characterization.

While this allows us to ignore complexes that are not color-shifted, there are still far too many such complexes for this to characterize the flag f-vectors of color-shifted complexes without a lot of additional work.  In the simple case of only two colors, if there are enough vertices, the number of color-shifted complexes with $n$ edges is the $n$-th partition number, which is greater than $13^{\sqrt{n}}$ for sufficiently large $n$.  Clearly, this gets too big for a brute force approach very quickly.

One may like to deal with this situation by putting additional restrictions on the complexes to consider.  This is what Frankl, F\"{u}redi, and Kalai \cite{ffk} did to characterize the f-vectors of colored complexes.  They showed that if there is an $n$-colored complex with a given f-vector, then there is an $n$-colored ``rev-lex" complex with the same f-vector.  Since there is only one possible $n$-colored rev-lex complex with a given f-vector, it is possible to try to construct that one particular complex, see whether it is a valid simplicial complex (as determined by checking some inequalities), and have that solve the problem of whether a given prospective f-vector corresponds to some $n$-colored complex.  Their approach echoed that of the Kruskal-Katona theorem \cite{kruskal, katona}, which showed that in order to characterize the f-vectors of all simplicial complexes, it sufficed to consider only the ``rev-lex" complexes.

One might hope to characterize the flag f-vectors of colored complexes in the same manner.  The major result of this paper is to say that this cannot be done:  we cannot put additional restrictions on what an extremal complex can look like locally beyond requiring that it be color-shifted.  We need some definitions in order to state the result.

\begin{definition}
\textup{Let $\Delta$ be a color-shifted, $n$-colored simplicial complex. A face $F \in \Delta$ is \textit{shift-maximal} if $\Delta - F$ is also a color-shifted, $n$-colored complex.}
\end{definition}

Note that this definition puts two conditions on a face in order for it be shift-maximal.  The face must be maximal with respect to inclusion, so that $\Delta - F$ is a simplicial complex.  It must also be maximal with respect to color-shifting, so that $\Delta - F$ is color-shifted.

\begin{definition}
\textup{Let $\Delta$ be an $n$-colored simplicial complex and let $S \subset [n]$.  The \textit{color-selected subcomplex} of $\Delta$ for the color set $S$ is the simplicial complex whose faces are precisely the faces of $\Delta$ whose color set is a subset of $S$.}
\end{definition}

By keeping the same coloring as $\Delta$, the color-selected subcomplex is an $|S|$-colored simplicial complex.

\begin{theorem}  \label{nolocal}
Let $\Delta$ be a color-shifted, $n$-colored simplicial complex.  There is an $m$-colored simplicial complex $\Gamma$ such that
\begin{enumerate}
\item the color-selected subcomplex of $\Gamma$ for the color set [n] is $\Delta$, and
\item $\Gamma$ is the unique color-shifted, $m$-colored simplicial complex with the flag f-vector $f(\Gamma)$.
\end{enumerate}
\end{theorem}

\proof  Let the shift-maximal faces of $\Delta$ be $F_1, F_2, \dots , F_k$.  Let $m = n+k$.  Let $\Gamma_i \subset \Delta$ be the complex whose unique shift-maximal face is $F_i$.  Let $w_i$ be the first vertex of color $n+i$ for $1 \leq i \leq k$.  Define $\Gamma$ as the union of $\Delta$ with a cone over $\Gamma_i$ with apex vertex $w_i$ for each $1 \leq i \leq k$.  Since all faces of $\Gamma$ not in $\Delta$ have one of the vertices $w_i$, which is of color $n+i \not \in [n]$, it is clear that $\Gamma$ satisfies the first property of the theorem.

Let $\Sigma$ be a color-shifted, $m$-colored simplicial complex such that $f(\Sigma) = f(\Gamma)$.  It suffices to show that $\Sigma = \Gamma$, as this will prove the second property of $\Gamma$.

Let $F_p$ be one of shift-maximal faces of $\Delta$, and let its vertex set be $F_p = \{v_{a_1}^{i_1}, v_{a_2}^{i_2}, \dots , v_{a_j}^{i_j}\}$.  We can compute that for any $r \in \{i_1, i_2, \dots , i_j\}$, $f_{\{r,p+n\}}(\Gamma) = a_r$.  Additionally, we can compute $f_{\{p+n\}}(\Gamma) = 1$.

Since $f_{\{p+n\}}(\Sigma) = f_{\{p+n\}}(\Gamma) = 1$, $\Sigma$ has only one vertex of color $p+n$.  For any $r \in \{i_1, i_2, \dots , i_j\}$, $f_{\{r, p+n\}}(\Sigma) = a_r$.  Since $\Sigma$ is color-shifted and there is only one way to arrange a given number of color-shifted edges on a set of two colors with only one vertex of one of the colors, $\Sigma$ has the same edges containing color $p+n$ as $\Gamma$.  Observe that $\Gamma$ has every possible face on the color set $\{i_1, i_2, \dots , i_j, p+n \}$ without requiring an edge that the complex is known not to have.  Because every edge of $\Sigma$ containing $w_p$ is an edge of $\Gamma$, in order for $f_{\{i_1, i_2, \dots , i_j, p+n \}}(\Sigma) = f_{\{i_1, i_2, \dots , i_j, p+n \}}(\Gamma)$, $\Sigma$ must likewise have all of these faces.  In particular, it follows from this that $F_p \in \Sigma$.  Furthermore, every face of $\Gamma$ containing a vertex of color $p+n$ is a subset of one of these faces.  Since $\Sigma$ is a simplicial complex, it must have all of these faces as well.

One can repeat this for all other $p \in [k]$ to show that the faces of $\Sigma$ containing a color outside of $[n]$ include all of the faces of $\Gamma$ containing a color outside of $[n]$.  Furthermore, $F_p \in \Sigma$ for all $p \in [k]$.  Because $\Sigma$ is color-shifted, this requires $\Gamma_p \subset \Sigma$ for all $p \in [k]$.  For every face $F \in \Gamma$ whose color set is contained in $[n]$, $F \in \Gamma_p$ for some $p \in [k]$.  Therefore, $F \in \Gamma_p \subset \Sigma$, and so $F \in \Sigma$.  We have shown that every face of $\Gamma$ is also a face of $\Sigma$, so $\Gamma \subset \Sigma$.  Since $f(\Sigma) = f(\Gamma)$, $\Sigma$ cannot have any other faces, and so $\Sigma = \Gamma$. \endproof

This result perhaps needs some explanation.  The result of Bj\"{o}rner, Frankl, and Stanley \cite{bfs} said that if we want to characterize the flag f-vectors of colored complexes, we can discard all complexes that are not color-shifted, and still have at least one color-shifted complex for each flag f-vector that has a corresponding simplicial complex.

Being color-shifted is essentially a local property in the sense that in order for a simplicial complex not to be color-shifted, the faces of some particular color set $S$ must not be color shifted.  The result of Bj\"{o}rner, Frankl, and Stanley says that rather than considering all of the ways to arrange the faces of color set $S$, we only need to consider the ones that satisfy the color-shifting condition.

Theorem~\ref{nolocal} says that if we exclude the ways to arrange faces of color set $S$ that are not color-shifted, then we cannot exclude any other ways to arrange the faces of color set $S$.  Suppose that we pick a color set $S$ and wish to exclude all complexes that are not color-shifted.  Suppose that we also wish to eliminate one particular way to arrange the faces of color set $S$ that is color-shifted.  We can fill in the lower dimensional faces whose color sets are subsets of $S$ in any color-shifted manner that we like, and the faces of color sets that are subsets of $S$ form a colored complex $\Delta$ with color set $S$.

Theorem~\ref{nolocal} says that there is a colored complex $\Gamma$ that meets the conditions of the theorem.  Hence, the color-selected subcomplex of $\Gamma$ with color set $S$ is $\Delta$, so $\Gamma$ has exactly the arrangement of faces of color set $S$ that we wished to exclude.  Restricting to the set of colored complexes that are color-shifted and whose faces of color set $S$ are not arranged in exactly the manner of $\Delta$ thus excludes $\Gamma$.  Because $\Gamma$ is the unique color-shifted complex with flag f-vector $f(\Gamma)$, this means that our slightly narrowed class of complexes has no complex with flag f-vector $f(\Gamma)$.  However, $\Gamma$ has flag f-vector $f(\Gamma)$, so our new class of complexes has a different characterization of the flag f-vectors than that of the color-shifted complexes.  The effort to restrict the class of complexes by eliminating the one particular arrangement of faces of color set $S$ without changing the characterization of the flag f-vectors is a failure.

Intuitively, this means that there is no hope of making progress toward a solution to this problem by putting stronger local restrictions on the class of complexes that must be considered, in addition to restricting to color-shifted complexes.  If there are to be additional local restrictions on the class of complexes that we consider, then we have to give up the color-shifting property.

\end{document}